\documentclass[11pt,leqno]{amsart}
\usepackage[T1]{fontenc}
\usepackage[utf8]{inputenc}
\usepackage[english]{babel}

\usepackage{amsmath, amsthm, amssymb}
\usepackage{esint} 

\usepackage{dsfont}

\usepackage{hyperref}
\usepackage{theoremref}
\usepackage{graphicx}
\usepackage{enumerate}

\usepackage{csquotes} 
\usepackage[
backend=biber,
style=alphabetic,
maxbibnames=99,
maxalphanames=5
]{biblatex}
\usepackage{thmtools}
\theoremstyle{plain}
\declaretheorem[title=Theorem, parent=section]{thm}
\declaretheorem[title=Lemma,sibling=thm]{lem}

\theoremstyle{definition}

\declaretheorem[title=Remark, numbered=no]{remark*}

\numberwithin{equation}{section}

\usepackage[backgroundcolor=white, bordercolor=blue,linecolor=blue]{todonotes}

\newcommand{\R}{\mathbb{R}}

\newcommand{\eps}{\varepsilon}
\newcommand{\loc}{\mathrm{loc}}

\DeclareMathOperator{\dist}{dist}

\DeclareMathOperator{\supp}{supp}

\newcommand{\dd}{\mathop{}\!\mathrm d} 

\begin{document}

\allowdisplaybreaks
\title{A note on a nonlocal antisymmetric barrier}
 
\author{Serena Dipierro}
\address{Department of Mathematics and Statistics,
University of Western Australia,
35 Stirling Highway, Crawley WA 6009, Australia}
\email{serena.dipierro@uwa.edu.au}

\author{Giorgio Poggesi}
\address{School of Mathematical Sciences, Adelaide University, Adelaide, SA 5005, Australia}
\email{giorgio.poggesi@adelaide.edu.au}

\author{Jack Thompson}

\address{Department of Mathematics, KTH, Lindstedsv\"agen 25, 100 44 Stockholm, Sweden}
\email{jackth@kth.se}

\author{Enrico Valdinoci}
\address{Department of Mathematics and Statistics,
University of Western Australia,
35 Stirling Highway, Crawley WA 6009, Australia}
\email{enrico.valdinoci@uwa.edu.au}



\allowdisplaybreaks

\begin{abstract}
In this note, we present a barrier function for an equation involving the fractional Laplacian in the antisymmetric setting, providing an explicit inequality in the full positivity set
of the function. This leads to an improvement of the quantitative maximum principle given in~\cite{DPTV_quantitative_2026},
where the bound was obtained in a strict subset of its positivity set and required additional assumptions on lower order coefficients. 

As a consequence of this, we give
an optimal lower bound on the growth of antisymmetric
solutions from the symmetry hyperplane.
\end{abstract}

\allowdisplaybreaks 
\maketitle

\section{Introduction}

Recently, in \cite{DPTV_quantitative_2026}, the current authors established a stability estimate for the fractional Serrin overdetermined problem. The proof relied on establishing several quantitative maximum principles for fractional equations with solutions that are \emph{antisymmetric}, namely, they are odd functions with respect to reflections across a given plane. In turn, this relied on the construction of a new antisymmetric barrier satisfying a particular upper bound for its fractional Laplacian, see \cite[Theorem 1.4]{DPTV_quantitative_2026}. In \cite{DPTV_quantitative_2026}, we proved that, in its positivity set, the barrier satisfies this upper bound pointwise away from a set of measure zero and it satisfies the bound distributionally {\em in a strict subset of its positivity set}. However, in our current research, it is useful to establish that
the barrier satisfies the estimate distributionally {\em in its full positivity set}. The purpose of this note is to prove this result. The result is as follows\footnote{In the statement of \cite[Theorem 1.4]{DPTV_quantitative_2026} there is a typo---Equation (1.11) should read \((-\Delta)^s \varphi \leq \frac{2(n+2s)}n x_1 \) in \(B_\rho^+(a)\setminus \partial B_\rho(a_\ast)\). This does not affect the results of \cite{DPTV_quantitative_2026}. We use the corrected constant throughout this note. \label{footnote:correction} }. 

\begin{thm} \thlabel{thm:main}
Let \(s\in (0,1)\), \(n\) be a positive integer, \(a=(a_1,\dots,a_n)\in \R^n\) with \(a_1\geq 0\), and \(\rho>0\). For \(x_0\in \R^n\), consider \begin{align*}
	\psi_{B_\rho(x_0)} (x) = \gamma_{n,s} \big (\rho^2 - \vert x- x_0\vert^2 \big)_+^s \qquad \text{where}\qquad \gamma_{n,s} = \frac{4^{-s} \Gamma \big(  \frac n2 \big ) }{\Gamma \big(  \frac {n+2s}2 \big )\Gamma ( 1+s )}
\end{align*} and set \begin{align*}
	\varphi(x) = x_1 \big ( \psi_{B_\rho(a)} (x)+ \psi_{B_\rho(a_\ast)} (x) )
\end{align*} where \(a_\ast =(-a_1,a_2,\dots,a_n) \) is the reflection of \(a\) across \(\{x_1=0\}\). Then \begin{align*}
	(-\Delta)^s \varphi \leq \frac{2(n+2s)}n x_1 
\end{align*} in \(B^+_\rho(a)\) in the distributional sense. 
\end{thm}

Here and throughout the paper, \(B_\rho(a)\) denotes the open ball in \(\R^n\) centred at \(a\) with radius \(\rho>0\), and 
$$B^+_\rho(a):=B_\rho(a)\cap\{x_1>0\}.$$ Furthermore, recall that, given an open set \(\Omega \subset \R^n\), \(u\in H^s_{\loc}(\Omega)\cap L^\infty(\R^n)\), and \(f\in L^1_{\loc}(\Omega)\), we say that \((-\Delta)^s u \leq f \) in \(\Omega\) in the distributional sense if \begin{align*}
	\frac{c_{n,s}} 2 \iint_{\R^n\times \R^n} \frac{(u(x)-u(y))(v(x)-v(y))}{\vert x-y\vert^{n+2s}} \dd y \dd x \leq \int_\Omega f v \dd x 
\end{align*} for all \(v\in C^\infty_c(\Omega)\) with \(v \geq 0\). Here \(c_{n,s}>0\) is a normalisation constant given by \begin{align}
	c_{n,s} = \frac{4^s s \Gamma \big ( \frac{n+2s}2 \big )}{\pi^{n/2}\Gamma(1-s)}. \label{eq:constant}
\end{align}

\begin{figure}[ht]
\centering
\fbox{
    \begin{minipage}{0.85\textwidth}
        \centering
     \fbox{   \includegraphics[width=0.95\textwidth]{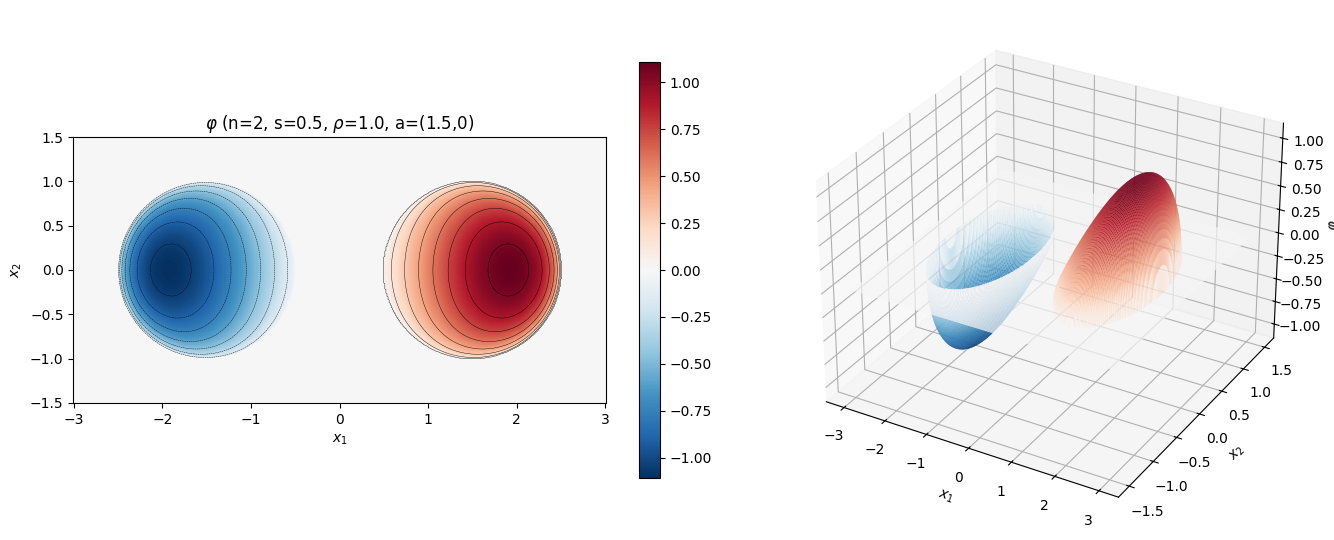} }\\
        \vspace{1em}
     \fbox{   \includegraphics[width=0.95\textwidth]{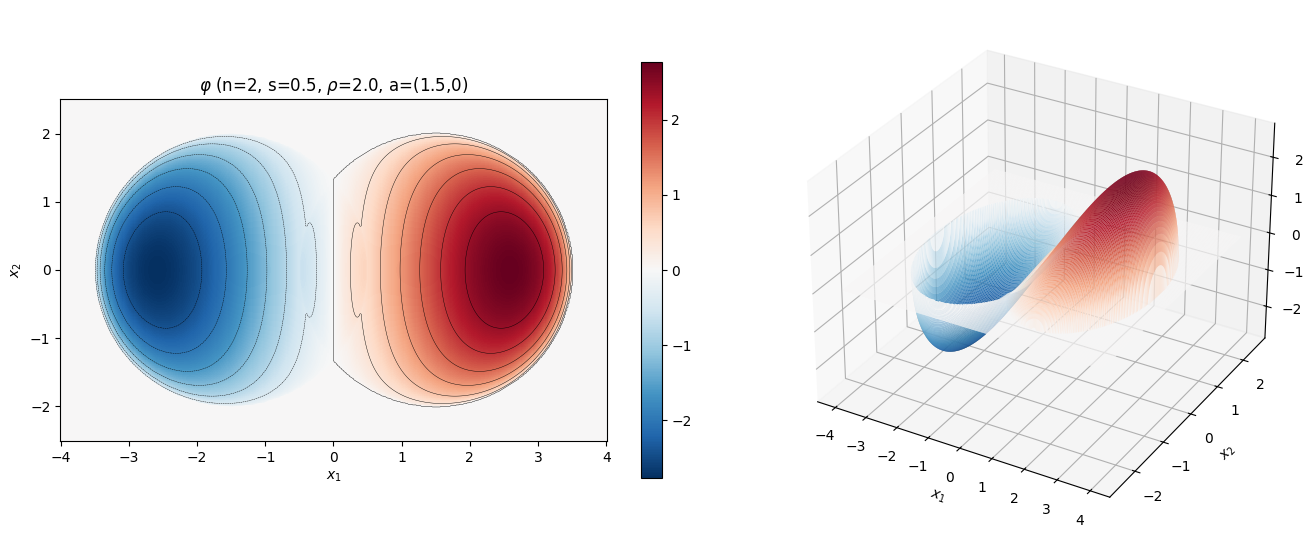} }
    \end{minipage}
}
\caption{Examples of the function~$\varphi$ in Theorem~\ref{thm:main}.}\label{GFif}
\end{figure}

The proof of Theorem~\ref{thm:main} will be presented in Section~\ref{THEPROOF}. See Figure~\ref{GFif} for a visualization of the function analyzed in Theorem~\ref{thm:main}.

In Section~\ref{Improved Quantitative MP:S}, we use Theorem~\ref{thm:main} to prove a generalisation of the quantitative maximum principle for antisymmetric functions established in \cite[Proposition 2.4]{DPTV_quantitative_2026}, given precisely in Theorem~\ref{Improved Quantitative MP}. The theorem allows us to obtain a lower bound for antisymmetric solutions
of equations involving the fractional Laplacian. Results of this kind are important since they detect the optimal growth from the symmetry hyperplane of the solution and enable the use of moving plane methods (see~\cite{DPTV_quantitative_2026}). 

Theorem~\ref{Improved Quantitative MP} generalises \cite[Proposition 2.4]{DPTV_quantitative_2026} in two ways. First, for technical reasons, the inequality obtained in \cite[Proposition 2.4]{DPTV_quantitative_2026} does not hold in a strip parallel to the plane of symmetry. In Theorem~\ref{Improved Quantitative MP}, by utilising Theorem~\ref{thm:main}, we are able to remove this artificial restriction, giving the inequality in the full half-space. Second, by a more careful analysis, we remove a smallness assumption on the zero-th order term present in \cite[Proposition 2.4]{DPTV_quantitative_2026}. Finally, we remark that \thref{thm:main} not only allows us to generalise \cite[Proposition 2.4]{DPTV_quantitative_2026}, but it greatly simplifies the proof, allowing us to follow are more concise and unified strategy than the original proof presented in \cite[Proposition 2.4]{DPTV_quantitative_2026}.

\subsection*{Acknowlegdements}

SD has been supported by the Australian Research Council (ARC) Future Fellowship FT230100333.
GP has been supported by the Australian Research Council (ARC) Discovery Early Career Researcher Award (DECRA) DE230100954.
EV has been supported by the Australian Research Council (ARC) Laureate Fellowship FL190100081.

\section{Proof of Theorem~\ref{thm:main}}\label{THEPROOF}

In our setting, it is convenient to first calculate an ``approximate fractional
Laplacian'' by removing a small ball, say of radius~$\eps$, centered at the singularity of the kernel, establish estimates uniform in \(\eps\), and then pass to the limit.
To this end,
given \(u\in L^\infty(\R^n)\) and \(\eps>0\), we define \begin{align*}
	(-\Delta)^s_\eps u(x) = c_{n,s}\int_{\R^n\setminus B_\eps (x)} \frac{u(x)-u(y)}{\vert x-y\vert^{n+2s}} \dd y  
\end{align*} where \(c_{n,s}\) is given by~\eqref{eq:constant}. By symmetry, we have \begin{align}
(-\Delta)^s_\eps u(x) = \frac{c_{n,s}}2 \int_{\R^n\setminus B_\eps } \big (2u(x)-u(x+y)-u(x-y) \big ) \frac{\dd y }{\vert y\vert^{n+2s}} . \label{eq:eps-frac-Lap}
\end{align}
In this framework, we have:

\begin{lem} \thlabel{thm:frac-Lap-phi-est}
	Let \(\varphi\) be as in \thref{thm:main}. Then, for all \(K\subset \subset B^+_\rho(a) \), \begin{align*}
		\sup_{\eps>0 } \vert (-\Delta)^s_\eps \varphi(x) \vert \leq C \big (1+\dist(x,\partial B_\rho(a_\ast))^{-s} \big ) \qquad \text{for all }x\in K\setminus \partial B_\rho(a_\ast).
	\end{align*}The constant \(C>0\) depends only on \(n\), \(s\), \(\rho\), \(a\), and \(K\).
\end{lem} 

\begin{proof}
Throughout the proof, \(C\) will always denote a positive constant depending only on \(n\), \(s\), \(\rho\), \(a\), and \(K\) whose value may change from line to line. Let \(\varphi = \varphi_1+\varphi_2\) where \(\varphi_1 = x_1 \psi_{B_\rho(a)}\) and \(\varphi_2 = x_1 \psi_{B_\rho(a_\ast)}\). Since \(\varphi_1\in C^\infty(B_\rho(a)) \cap L^\infty(\R^n)\), by standard estimates, we have that \begin{align*}
	\sup_{\eps>0 }\vert (-\Delta)^s_\eps \varphi_1(x) \vert \leq C \qquad \text{for all }x\in K.
\end{align*} Indeed, using~\eqref{eq:eps-frac-Lap}, if \(K'\) is such that \(K\subset \subset K'\subset \subset B_\rho(a)\) and \(r=\dist (x,\partial K')\) then  \begin{align*}
	\vert (-\Delta)^s_\eps \varphi_1(x) \vert &\leq C \bigg (  \| D^2\varphi_1\|_{L^\infty (K') }\int_{B_r} \frac{\dd y }{\vert y \vert^{n+2s-2}}+ \| \varphi_1 \|_{L^\infty(\R^n)}\int_{\R^n\setminus B_r} \frac{\dd y }{\vert y \vert^{n+2s}} \bigg ) \\
	&\leq C \big ( r^{2(1-s)}  + r^{-2s} ) \\
	&\leq C 
\end{align*} where the last inequality follows from the estimate \( \dist (K,\partial K') \leq r\leq \rho  \).

Next, we estimate \((-\Delta)^s_\eps\varphi_2\). Let \(\delta(x) =\dist(x,\partial B_\rho(a_\ast) )\). By a direct computation, we have that \begin{align}
	\vert  D^2 \varphi_2(x) \vert \leq C \delta(x)^{s-2} \qquad \text{for all }x\in \R^n\setminus \partial B_\rho(a_\ast). \label{eq:phi2-est}
\end{align} Now, if \(x\in K \setminus \partial B_\rho(a_\ast)\), we write \(r(x) = \min \big \{ \frac{\delta(x)} 2 , 1 \big  \}\) and \begin{align*}
	\vert (-\Delta)^s_\eps\varphi_2(x)\vert  &\leq C \int_{ \R^n  } \big \vert 2\varphi_2(x)-\varphi_2(x+y)-\varphi_2 (x-y) \big \vert  \frac{\dd y }{\vert y\vert^{n+2s}} \\
	&=C( I_1 + I_2 + I_3 ) 
\end{align*} where \begin{align*}
	I_1 &= \int_{B_{r(x)/2} } \big \vert 2\varphi_2(x)-\varphi_2(x+y)-\varphi_2(x-y) \big \vert  \frac{\dd y }{\vert y\vert^{n+2s}}\\
	I_2 &=  \int_{B_1 \setminus B_{r/2}  } \big \vert 2\varphi_2(x)-\varphi_2(x+y)-\varphi_2(x-y) \big \vert \frac{\dd y }{\vert y\vert^{n+2s}}; \text{ and}\\
	I_3 &=  \int_{\R^n \setminus B_1 } \big \vert 2\varphi_2(x)-\varphi_2(x+y)-\varphi_2(x-y) \big \vert \frac{\dd y }{\vert y\vert^{n+2s}}
\end{align*} For \(I_1\), we use~\eqref{eq:phi2-est} to obtain \begin{align*}
	I_1 &\leq C \delta(x)^{s-2}  \int_{ B_{r(x)/2}  } \frac{\dd y }{\vert y \vert^{n+2s-2}} \leq  C \delta(x)^{s-2}  \int_0^{r(x)/2 } \frac{\dd t }{t^{2s-1}} \leq C \delta(x)^{s-2}r(x)^{2-2s} \leq C \delta(x)^{-s}.
\end{align*} For \(I_2\), we have that \(\varphi_2\in C^s(\R^n)\), so \begin{align*}
	I_2 &\leq C [\varphi_2 ]_{C^s(\R^n)} \int_{B_1 \setminus B_{r(x)/2}  } \frac{\dd y }{\vert y\vert^{n+s}} \leq C \int^1_{r(x)/2}  \frac{\dd t }{t^{1+s}} \leq C r(x)^{-s}\leq C (1+\delta (x)^{-s}).
\end{align*} For \(I_3\), we have that \begin{align*}
	I_3 \leq C \| \varphi_2\|_{L^\infty(\R^n)} \int_{\R^n \setminus B_1 } \frac{\dd y }{\vert y\vert^{n+2s}} \leq C . 
\end{align*} Collecting the above inequalities, we obtain \begin{align*}
	\sup_{\eps>0 }\vert (-\Delta)^s_\eps \varphi_2(x) \vert \leq C \big ( 1+ \delta(x)^{-s} ) \qquad \text{for all }x\in K\setminus \partial B_\rho(a_\ast). 
\end{align*} Recalling the definition of \(\delta(x)\) and using that \begin{align*}
	\sup_{\eps>0 }\vert (-\Delta)^s_\eps \varphi(x) \vert \leq \sup_{\eps>0 }\vert (-\Delta)^s_\eps \varphi_1(x) \vert +\sup_{\eps>0 }\vert (-\Delta)^s_\eps \varphi_2(x) \vert,
\end{align*} the proof is complete. 
\end{proof}

We now give the proof of \thref{thm:main}. 

\begin{proof}[Proof of \thref{thm:main}]
	Let \(v\in C^\infty_c(B_\rho^+(a))\) with \(v\geq 0\).  On one hand, since \(\varphi, v \in H^s(\R^n)\), we have that \begin{align*}
		\frac{c_{n,s}} 2 \iint_{\R^n\times \R^n } &\frac{(\varphi(x)-\varphi(y))(v(x)-v(y))}{\vert x-y\vert^{n+2s}} \dd y \dd x \\
		&= \frac{c_{n,s}} 2  \lim_{\eps \to 0 }\iint_{(\R^n \times \R^n ) \cap \{ \vert x-y\vert >\eps \}} \frac{(\varphi(x)-\varphi(y))(v(x)-v(y))}{\vert x-y\vert^{n+2s}} \dd y \dd x
	\end{align*} On the other hand, for all \(\eps>0\), we have that \begin{align*}
		\frac{c_{n,s}} 2 \int_{\R^{2n} \cap \{ \vert x-y\vert >\eps \}} \frac{(\varphi(x)-\varphi(y))(v(x)-v(y))}{\vert x-y\vert^{n+2s}} \dd y \dd x &= \int_{B_\rho^+(a)} v (-\Delta)^s_\eps \varphi \dd x . 
	\end{align*} Since the limit \begin{align*}
		(-\Delta)^s \varphi(x) &= \lim_{\eps \to 0} (-\Delta)^s_\eps  \varphi(x)
	\end{align*} exists for all \(x\in B_\rho^+ (a) \setminus \partial B_\rho(a_\ast )\), the set \(B_\rho^+ (a) \cap  \partial B_\rho(a_\ast )\) has zero Lebesgue measure, and, by \thref{thm:frac-Lap-phi-est} with \(K=\supp v\), we have \begin{align*}
		\vert v (-\Delta)^s_\eps \varphi \vert \leq C \big (1+\dist(\cdot ,\partial B_\rho(a_\ast))^{-s} \big ) \in L^1 (B_\rho^+(a)),
	\end{align*} the dominated convergence theorem implies that \begin{align*}
		\lim_{\eps \to 0}\int_{B_\rho^+(a)} v (-\Delta)^s_\eps \varphi \dd x  = \int_{B_\rho^+(a)} v (-\Delta)^s \varphi \dd x . 
	\end{align*} Hence, \begin{align*}
		\frac{c_{n,s}} 2 \iint_{\R^n\times \R^n } \frac{(\varphi(x)-\varphi(y))(v(x)-v(y))}{\vert x-y\vert^{n+2s}} \dd y \dd x &=  \int_{B_\rho^+(a)} v (-\Delta)^s \varphi \dd x.  
	\end{align*} Finally, in \cite[Theorem 1.4]{DPTV_quantitative_2026} (see also Footnote~\ref{footnote:correction}), it was proved that \begin{align*}
		(-\Delta)^s \varphi (x) \leq \frac{2(n+2s)}nx_1 \qquad \text{for all } x\in B_\rho^+ (a) \setminus \partial B_\rho(a_\ast ), 
	\end{align*} provided \(a_1>0\). This inequality also remains valid for \(a_1=0\) since, in this case, \(\varphi(x) = 2x_1\psi_{B_\rho}(a) \), so \((-\Delta)^s \varphi(x) = \frac{2(n+2s)}nx_1 \) in \(B^+_\rho\) by \cite[Lemma 2.2]{DPTV_quantitative_2026}. Hence, we obtain \begin{align*}
		\frac{c_{n,s}} 2 \iint_{\R^n\times \R^n } \frac{(\varphi(x)-\varphi(y))(v(x)-v(y))}{\vert x-y\vert^{n+2s}} \dd y \dd x &\leq  \frac{2(n+2s)}n \int_{B_\rho^+(a)} x_1 v  \dd x
	\end{align*} as required. 
\end{proof}

\section{An application}\label{Improved Quantitative MP:S}

We now provide an application of \thref{thm:main}: we show that \thref{thm:main} can be used to give a more direct and concise proof of the quantitative maximum principle for antisymmetric functions established in \cite[Proposition 2.4]{DPTV_quantitative_2026}. In fact, the proof in \cite{DPTV_quantitative_2026} required the analysis of three separate cases, whereas the proof we provide here allows a unified analysis, thanks to \thref{thm:main}. Not only is the new proof more concise, but it also allows us to establish a more general version\footnote{Specifically,
the conclusion of~\cite[Proposition 2.4]{DPTV_quantitative_2026} consisted in the inequality~\eqref{eq:quantitative maximum principle}, but with validity restricted
to the cases in which the distance of~$a$ to the boundary of~$H$
lied in~$\left[0,\frac\rho2\right]\cup(\rho,+\infty)$. Theorem~\ref{Improved Quantitative MP} removes this artificial restriction, giving
the full validity of the lower bound in~\eqref{eq:quantitative maximum principle}.} of \cite[Proposition 2.4]{DPTV_quantitative_2026}. Furthermore, by a more careful analysis, we only require that the zero-th order term in Equation~\eqref{PERv} be bounded, that is, we do not require the smallness assumption (2.7) in \cite[Proposition 2.4]{DPTV_quantitative_2026}. We present the result in \thref{Improved Quantitative MP} below.

For the reader's convenience, we adopt some notation from \cite{DPTV_quantitative_2026}.
Given an open half-space $H\subset \mathbb{R}^n$, we denote by $Q_{\partial H} : \mathbb{R}^n \to \mathbb{R}^n$ the function that reflects $x$ across the hyperplane $\partial H$, and we say that a function $v: \mathbb{R}^n \to \mathbb{R}$ is antisymmetric with respect to $\partial H$ if $v(x) = - v \left( Q_{\partial H} \left( x \right) \right)$ for all $x \in \mathbb{R}^n$.

\begin{thm}\thlabel{Improved Quantitative MP}
	Let $ H \subset \R^n $ be an open half-space and $U $ be an open subset of $H$.
	 Let~$c \in L^\infty(U)$. 
	
	 Let $a\in \overline{H} $ and $ \rho>0 $ be such that $B_\rho(a)\cap H \subset U $, and set $ a_\ast := Q_{\partial H}(a) $ be the reflection of~$a$ across~$\partial H$. 
	
	Let $K\subset H$ be a non-empty bounded open set that is disjoint from $B_\rho(a)$ and let $0 \le \mathcal{S} <\infty$ be such that
	\begin{equation}\label{eq:Sup def}
		\sup_{\genfrac{}{}{0pt}{2}{x\in K}{y\in B_\rho(a)\cap H}}
		\vert Q_{\partial H}(x)-y\vert \leq \mathcal{S} .
	\end{equation}
	
	If~$v\in H^s(\R^n)$ is antisymmetric with respect to~$\partial H$ and satisfies
	\begin{equation}\label{PERv}
		\begin{cases}
			(-\Delta)^s v +c v \geq 0 \quad & \text{ in } U  ,\\
			v \geq 0 \quad & \text{ in } H ,
		\end{cases}
	\end{equation}
	in the weak sense, then
	we have that
	\begin{equation}\label{eq:quantitative maximum principle}
		v(x) \geq C \| \delta_{\partial H} v\|_{L^1(K)} \, \delta_{\partial H}(x) \, \left( \rho^2 - |x-a|^2 \right)_+^s \quad \text{ for a.e. } x \in H \cap  B_{\rho}(a) ,
	\end{equation}
	with
	\begin{equation*}
		C := \frac{ C(n,s)}{\mathcal{S}^{n+2s+2}} \Big(1+\rho^{2s}\|c^+\|_{L^\infty(U)} \Big)^{-1}.
	\end{equation*}
\end{thm}
\begin{proof}
	Without loss of generality, we can assume that
	$$
	H = \R^n_+ = \big\{ x=(x_1,\dots,x_n)\in\R^n \, \text{ : } \, x_1>0 \big\} ,
	$$
	and use the notations
	\begin{equation*}
	A^+ = A \cap \R^n_+ , \text{ for all } A\subset \R^n , 
	\end{equation*}
	and
	\begin{equation*}
	x_\ast:= Q_{\partial H }(x) = Q_{\partial \R^n_+ }(x)=x-2x_1e_1 , \text{ for all } x\in\R^n .
	\end{equation*}
	\smallskip
	
	Set $B:=B_\rho(a)$ and $B_\ast:=B_\rho( a_\ast )$.
	Let $\tau \geq 0$ be a constant to be chosen later and set
	$$
	w := \tau \varphi +  (\mathds{1}_K + \mathds{1}_{K_\ast}) v ,
	$$
	where $K_\ast:=Q_{\partial \R^n_+}(K)$ and $\varphi$ is as defined in \thref{thm:main}. 
	
	To comply with the weak form of the equation,
	we perform integral calculations against smooth test functions compactly supported in $B^+$
	 and take a smooth approximations when necessary. 
	  By \thref{thm:main}, for any $\xi \in C^\infty_0(B^+ )$ with $\xi \geqslant 0$, we have that 
	\begin{equation*}
		\begin{split}
			\mathcal E (w, \xi ) &:=
			\frac{c_{n,s}}2
\iint_{\R^n\times \R^n} \frac{(w(x) - w(y))(\xi(x) - \xi(y))}{|x - y|^{n+2s}}\, dx\, dy
\\&			=  \tau \mathcal E (\varphi, \xi )
			+ \mathcal E (\mathds{1}_K v , \xi )
			+ \mathcal E (\mathds{1}_{K_\ast} v , \xi ) 
			\\
			& \leq  \frac{2(n+2s)}{n} \,  \tau \int\limits_{ \supp(\xi) } x_1 \xi(x) \dd x 
			- c_{n,s} \left[ \ \, 
			\iint\limits_{ \supp(\xi)\times K} \frac{ v(y) \xi(x) }{\vert x-y\vert^{n+2s}} \dd x  \dd y
			+
			\iint\limits_{ \supp(\xi)\times K_\ast} \frac{ v(y) \xi(x) }{\vert x-y\vert^{n+2s}} \dd x  \dd y\right]
			\\
			& = \int\limits_{\supp(\xi) } \bigg [ \frac{2(n+2s)}{n} \, \tau x_1 - c_{n,s}  \int\limits_{K}  \bigg ( \frac 1{\vert x-y\vert^{n+2s}} - \frac 1 {\vert x_\ast - y \vert^{n+2s}} \bigg ) v(y) \dd y \bigg ] \xi(x) \dd x ,
		\end{split}
	\end{equation*}
	where $ \supp(\xi) $ denotes the support of $\xi$.
	Note that 
	\begin{equation*}
		\bigg \vert \iint\limits_{ \supp(\xi)\times K} \frac{ v(y) \xi(x)}{\vert x-y\vert^{n+2s}} \dd x  \dd y\bigg \vert <+\infty .
	\end{equation*} 
	Indeed, since $B$ and $ K $ are disjoint and $\supp \xi$ is compactly contained in $B$, we see that $ \dist (\supp \xi , K)>~0 $, and hence 
	\begin{equation*}
		\begin{split}
			\bigg  \vert  \iint\limits_{ \supp(\xi)\times K} \frac{ v(y) \xi(x)}{\vert x-y\vert^{n+2s}} \dd x  \dd y \bigg \vert &\leq \dist (\supp \xi , K)^{-n-2s} \| \xi\|_{L^1(\R^n)} \| v\|_{L^1(K)} \\
			&\leq \dist (\supp \xi , K)^{-n-2s} \vert K \vert^{1/2} \| \xi\|_{L^1(\R^n)} \| v\|_{L^2(\R^n)}. 
		\end{split}
	\end{equation*} 
	A similar estimate can be obtained for the integral involving \(K_\ast\).
	
	Recalling \eqref{eq:Sup def} and noting that, for all $ x\in B^+ $ and $ y\in K $,
	\begin{equation*}
		\begin{split}
			\frac 1{\vert x-y\vert^{n+2s}} - \frac 1 {\vert x_\ast- y \vert^{n+2s}}  
			& = \frac{n+2s}2 \int\limits_{\vert x- y \vert^2}^{\vert x_\ast- y \vert^2} t^{-\frac{n+2s+2}2} \dd t
			\\
			& \geq \frac{n+2s}2  \;\frac{ \vert x_\ast - y \vert^2 - \vert x- y \vert^2 }{ \vert x_\ast- y \vert^{n+2s+2}}
			\\
			&= 2(n+2s)  \frac{x_1y_1}{\vert x_\ast -y \vert^{n+2s+2} } 
			\\
			&\geq \frac{2(n+2s)}{\mathcal{S}^{n+2s+2}} x_1y_1 ,
		\end{split}
	\end{equation*}
	we obtain that
	\begin{equation*}
		\begin{split}
			\mathcal E(w,\xi) 
			& \leq \left(
			\frac{2(n+2s)}n\,\tau -\frac{2c_{n,s}(n+2s)}{\mathcal{S}^{n+2s+2}}\,\|y_1v\|_{L^1(K)} \right)  \int\limits_{\supp(\xi) }   x_1 \xi(x) \dd x 
			\\
			& =
			C \left( \tau - \widetilde C \,  \frac{ \| y_1 v \|_{L^1(K)}}{\mathcal{S}^{n+2s+2}} \right)  \int\limits_{\supp(\xi) }   x_1 \xi(x) \dd x ,
		\end{split}
	\end{equation*}
	where~\(C\) and~\(\widetilde C\) depend only on~\(n\) and~\(s\).
	
	Hence, using that
	\begin{equation*}
		w = \tau \varphi = \tau x_1( \psi_B +\psi_{B_\ast} ) \leq 2 \, \gamma_{n,s} \rho^{2s} \tau x_1 \quad \text{ in } B^+ \supset \supp(\xi) , 
	\end{equation*}
	we have that
	\begin{equation*}
		\mathcal E(w,\xi) + \int\limits_{\supp(\xi) } c^+(x) w(x) \xi(x) \dd x  
		\leq
		C  \left[ \tau \left( 1  + \rho^{2s}\|c^+\|_{L^\infty(U)} \right)  - \widetilde C \, \frac{  \| y_1 v \|_{L^1(K)} }{\mathcal{S}^{n+2s+2}} \right] \int\limits_{\supp(\xi)}   x_1 \xi(x) \dd x ,
	\end{equation*}
	where~$C$ and~$\widetilde C$ may have changed from the previous formula but still depend only on~$n$ and~$s$. 
	
	We thus get that
	\begin{equation*}
		(-\Delta)^s w + c^+ w \leq 0 \quad \text{ in } B^+ ,
	\end{equation*}
	provided that \begin{equation*}
		\tau \le \frac{ \widetilde C\, \| y_1 v \|_{L^1(K)}}{2\mathcal{S}^{n+2s+2}} \left( 1  + \rho^{2s}\|c^+\|_{L^\infty(U)} \right)^{-1}
		. 
	\end{equation*}
	
	We now note that
	\begin{equation}\label{deiwty0987654}
		w = \mathds{1}_K v \le v \quad {\mbox{ in }}\R^n_+\setminus B^+ .
	\end{equation}
	
	Hence, setting
	\begin{equation}\label{eq:def:tau}
		\tau :=  \widetilde C   \, \frac{  \| y_1 v \|_{L^1(K)}}{2\mathcal{S}^{n+2s+2}} \left( 1  + \rho^{2s}\|c^+\|_{L^\infty(U)} \right)^{-1} ,
	\end{equation} 
	we have that
	\begin{equation*}
		\begin{cases}
			(-\Delta)^s w + c^+ w \leq 0 \quad & \text{ in } B^+ ,
			\\
			v \ge w \quad & \text{ in } \R^n_+ \setminus B^+  .
		\end{cases}
	\end{equation*}
	
	Let now~$\eta_j\in C_c^\infty(K)$ with
	$0\leq\eta_j\leq1$ be such that~$\eta_j\to1$ a.e. in $K$. We define
	\[
	 m_j:=\int_K \eta_j(y)y_1v(y)\,\dd y,\qquad
	 \tau_j:=\widetilde C\,\frac{m_j}{2\mathcal S^{n+2s+2}}
	 \left(1+\rho^{2s}\|c^+\|_{L^\infty(U)}\right)^{-1},
	\]
	and 
	\[
	 w_j:=\tau_j\varphi+(\eta_j+\eta_j\circ Q_{\partial H})v\in H^s(\R^n).
	\]
	Notice that~$w_j$ is
	antisymmetric. Moreover,
	since~$0\leq\eta_j\leq1$ and~$v\geq0$ in $H$,
	we have that~$w_j\leq v$ in~$H\setminus B^+$.
By repeating the above calculation
	with $\eta_jv$ in place of $\mathds{1}_Kv$,
	$\tau_j$ in place of
	$\tau$, and~$m_j$ in place of $\|y_1v\|_{L^1(K)}$, we see that, in the weak sense,
	\[
	 (-\Delta)^s w_j+c^+w_j\leq0\qquad\text{in }B^+.
	\]

Additionally, since $v\geq0$ in $H$ and $c^+\geq c$,
	\[
	 (-\Delta)^s v+c^+v\geq(-\Delta)^s v+cv\geq0
	 \qquad\text{in }B^+.
	\]
	Therefore $z_j:=v-w_j$ is antisymmetric, nonnegative in
	$H\setminus B^+$, and satisfies
	\[
	 (-\Delta)^s z_j\geq -c^+z_j\qquad\text{in }B^+.
	\]
	The comparison principle in \cite[Proposition 3.1]{FallJarohs2015}
	(see also \cite[Proposition 2.3]{JarohsWeth2019}) applies, since the
	coefficient $-c^+$ is bounded above by $0<\lambda_1(B^+)$. Hence
	\[
	 v\geq w_j=\tau_j\varphi\qquad\text{a.e. in }B^+.
	\]
	We now pass to the limit as~$j\to+\infty$ via the Dominated Convergence Theorem. In this way, we see that~$m_j\to\|y_1v\|_{L^1(K)}$,
	$\tau_j\to\tau$, where $\tau$ is given by~\eqref{eq:def:tau} and
	\[
	 v\geq\tau\varphi\qquad\text{a.e. in }B^+.
	\]
	
{F}rom this, using that
	\begin{equation*}
		\varphi(x) = x_1 (\psi_B(x) + \psi_{B_\ast}(x))  \ge x_1 \psi_B(x) ,
	\end{equation*}
	we obtain
	\begin{equation*}
		v(x) \geq \tau \, x_1 \, \psi_B(x) \quad \text{ in } B^+ ,
	\end{equation*}
	where $\tau$ is as defined in \eqref{eq:def:tau}. This is the desired result.
\end{proof}

\printbibliography

\vfill
\end{document}